\documentclass[12pt]{article}
\usepackage{amssymb}
\usepackage{amsfonts}

\input{tcilatex}
\begin{document}

\begin{center}
\bigskip

\textbf{Note on islands in path-length sequences of binary trees}
\end{center}

\bigskip

\begin{center}
Santiago Cortes Reina,* Stephan Foldes,*

Yousof Mardoukhi,* Navin M. Singhi**

\bigskip 

*Tampere University of Technology

PL 553, 33101 Tampere, Finland

sf@tut.fi

\bigskip 

**Indian Institute of Technology Bombay

Mumbai 400 076, India

navin@math.iitb.ac.in

\bigskip 

September 2014

\bigskip

\textbf{Abstract}
\end{center}

\textit{An earlier characterization of topologically ordered (lexicographic)
path-length sequences of binary trees is reformulated in terms of an
integrality condition on a scaled Kraft sum of certain subsequences (full
segments, or islands). The scaled Kraft sum is seen to count the set of
ancestors at a certain level of a set of topologically consecutive leaves is
a binary tree.}

\textit{\bigskip }


\bigskip 

\bigskip 

\bigskip 

\textbf{1 Introduction and connections with the theory of islands}

\bigskip

In the construction of optimal average-length uniquely decipherable codes,
Kraft's and McMillan's theorems guarantee that only instantaneous (pre-fix
free) codes need to be considered, and these correspond to binary trees.
More precisely, we consider finite rooted binary trees with a "topological"
structure specifying for each non-leaf node its left son and right son.
Alternatively, this topological structure can be thought of as a labelling
of the edges of the tree in such a way that the two edges from each non-leaf
node to its sons are labelled $0$ and $1$ (for left and right). The $0-1$
sequences corresponding to the root-to-leaf paths of the tree \ then
determine an instantaneous code, and all instantaneous codes arise this way.
In fact Huffman's classical algorithm constructs the code by constructing
the topological tree. Ordering the set of codewords lexicographically
corresponds to a topological left-to-right enumeration of the root-to-leaf
paths, and writing down the sequence of lenghts produces the \textit{%
lexicographic (topological) length sequence, }that in [FS] was shown to
fully determine the code.

In [FS] lexicographic length sequences were characterized among all
sequences of \textit{positive} (meaning here non-negative) integers,
including $0$. This characterization is reformulated in the present note,
and a certain integer parameter appearing in the characterization is given a
combinatorial meaning. The characterization is established by induction on
the number of nodes of an associated (in general not binary) tree.

Consider any finite sequence $l_{1},....,l_{n}$ of $n\geq 1$ positive real
numbers. For any non-empty sub-interval $[i,j]\subseteq \lbrack 1,n]$
consider the \textit{neighborhood maximin} parameter $m[i,j]$ defined by%
\[
m[i,j]=\max_{[i,j]\subset \lbrack h,k]\subseteq \lbrack 1,n]}\min
(l_{h},...,l_{k})
\]%
where the maximum is taken over all subintervals $[h,k]$\ of $[1,n]$
properly containing $[i,j]$, and the maximin is $0$ if $[i,j]=[1,n]$.
Reformulating a definition in [FS], where the neighborhood maximin paeameter
was called the \textit{largest value near} the interval, the interval $[i,j]$
is called a \textit{full segment }(or \textit{island} in the terminology of
Cz\'{e}dli [C] and [FHRW], where further references are also given) if $%
m[i,j]<\min (l_{i},...,l_{j}).$ As two full segments are either disjoint or
comparable by inclusion (a fact apparently first noted and used by Gernot H%
\"{a}rtel [G]), the set of all full segments is the node set of a rooted
tree, called \textit{tree of full segments}, where $[1,n]$ is the root and
for any other full segment $I$ the father of $I$ is the smallest full
segment containing $I$ properly. (Similarly, under mild assumptions on the
generalized island systems studied in [FHRW], islands are the nodes of a 
\textit{tree of islands}.) \ For any $[i,j]\subseteq \lbrack 1,n]$ we also
write $K[i,j]$ for the \textit{partial Kraft sum} $2^{-l_{i}}+...+2^{-l_{j}}.
$ For $[i,j]=[1,n]$ this is just the Kraft sum of $l.$

\bigskip

\textbf{2. Characterisation of topological length sequences of binary trees}%
\bigskip

Using the tree of islands the following variant formulation of the
characterization of topologically ordered path-length sequences given in
[FS] can be proved:

\bigskip

\textbf{Theorem} (variant of Theorem 4 of [FS]) \textit{A finite non-empty
sequence} $l=(l_{1},...,l_{n})$ \textit{of positive integers with Kraft sum} 
$1$\textit{\ is the path-length sequence of a topological binary tree if and
only if for every full segment} $S=[i,j]\subseteq \lbrack 1,n],$\textit{\
the number }$2^{m(S)}K(S)$ \textit{is an integer, where }$m(S)$ \textit{and} 
$K(S)$ \textit{denote the neighborhood maximin parameter and the partial
Kraft sum of} $S$ \textit{. In that case,} $2^{m(S)}K(S)$ \textit{is the
number of ancestors at level} $m$ \textit{below the root of the leaves} $%
v_{i},...,v_{j}$ \textit{(in the topologically ordered enumeration }$%
v_{1},...,v_{n}$ o\textit{f all leaves).}

\bigskip

\textbf{Proof.} \ Suppose $l$ is the path-length sequence of a binary tree
with a left-right topology. For each level $m$ ancestor of a leaf with index
in $S$, its descendents form a binary tree. If there are $k$ of these trees,
then each term in $K(S)$ corresponds to exactly one term in the Kraft sum of
the path-length sequence of exactly one of the $k$ binary trees, and it is
in value equal to this latter term divided by $2^{m}$.

\bigskip

Conversely, we can proceed by induction on the number of nodes of the tree
of full segments of a given sequence. The basis of induction is trivial.
Suppose that for the sequence of positive integers $l,$ for every full
segment $S=[i,j]\subseteq \lbrack 1,n]$\textit{\ }the number\textit{\ }$%
2^{m(S)}K(S)$ is an integer. Take a minimal full segment $S=[i,j].$\ Let $%
l^{\prime }$ denote the segment obtained by replacing in $l$ the -
necessarily constant - subsequence $l_{i},...,l_{j}$ by a term equal to $%
m[i,j]$ repeated $2^{m(S)}K(S)$ times. Then the following map $f$ is a
bijection from the set of full segments of $l$ with $S$ removed, to the set
of all full segments of $l^{\prime }$:%
\[
f:[h,k]\mapsto \lbrack s,t] 
\]%
where, denoting $2^{m(S)}K(S)$ by $r$,

\begin{eqnarray*}
\lbrack s,t] &=&[h,k]\text{ if }k<i \\
\lbrack s,t] &=&[h-(j-i+1)+r,k-(j-i+1)+r]\text{ if }j<h \\
\lbrack s,t] &=&[h,k-(j-i+1)+r]\text{ otherwise.}
\end{eqnarray*}%
The map preserves partial Kraft sums and neighborhood maximin parameters and
it also establishes an isomorphism between the tree of full segments of $l$
with the leaf node $S$ clipped, and the tree of all full segments of $%
l^{\prime }.$ Based on inductive hypothesis, construct the binary tree for $%
l^{\prime }$ then append, to each of the $r$ leaves starting from the $i$'th
leaf, uniform binary trees of depth $l_{i}-m(S).$ \ \ \ \ \ \ \ \ \ \ \ \ \
\ \ \ \ \ \ \ \ \ \ \ \ \ \ \ \ \ \ \ \ \ \ \ \ \ \ \ \ \ \ \ \ \ \ \ \ \ \
\ \ \ \ \ \ \ \ \ \ \ \ \ \ \ \ \ \ \ \ \ \ \ \ \ \ \ \ $\square $

\bigskip

\textbf{References}

\bigskip

[C] \ G. Cz\'{e}dli, The number of rectangular islands by means of
distributive lattices. \textit{European J. Combinatorics} 30 (2009), 208-215

\bigskip 

[FS] \ S. Foldes, N.M. Singhi, On instantaneous codes. \textit{J.
Combinatorics Inf. Syst. Sci.} 31 (2006) 307-317

\bigskip 

[FHRW] S. Foldes, E.K. Horvath, S. Radeleczki, T Waldhauser, A general
framework for island systems. To appear in \textit{Acta Sci. Math. (Szeged)}%
. Manuscript on ArXiv, 2013

\bigskip 

[H] \ G. H\"{a}rtel, personal communication, Tampere University of
Technology (2007)

\bigskip 

\end{document}